\documentclass{article}
\usepackage{amsmath,amssymb,amsthm,enumitem,hyperref,url,xcolor,zref-clever}
\usepackage[style=alphabetic,url=false,isbn=false,giveninits=true,maxbibnames=9]{biblatex}
\renewbibmacro{in:}{\ifentrytype{article}{}{\printtext{\bibstring{in}\intitlepunct}}}
\usepackage{authblk}
\usepackage{todonotes}

\setlist[enumerate,1]{label=(\roman*)}
\numberwithin{equation}{section}

\newtheorem{theorem}{Theorem}[section]
\zcRefTypeSetup{theorem}{cap}

\AddToHook{env/con/begin}{%
\zcsetup{countertype={theorem=con}}}
\zcRefTypeSetup{con}{cap}
\newtheorem{lemma}[theorem]{Lemma}
\AddToHook{env/lemma/begin}{%
\zcsetup{countertype={theorem=lemma}}}
\zcRefTypeSetup{lemma}{cap}
\newtheorem{proposition}[theorem]{Proposition}
\AddToHook{env/proposition/begin}{%
\zcsetup{countertype={theorem=proposition}}}
\zcRefTypeSetup{proposition}{cap}
\newtheorem{corollary}[theorem]{Corollary}
\AddToHook{env/corollary/begin}{%
\zcsetup{countertype={theorem=corollary}}}
\zcRefTypeSetup{corollary}{cap}

\AddToHook{env/observation/begin}{%
\zcsetup{countertype={theorem=observation}}}
\zcRefTypeSetup{observation}{cap}
\theoremstyle{definition}
\newtheorem{definition}[theorem]{Definition}
\AddToHook{env/definition/begin}{%
\zcsetup{countertype={theorem=definition}}}
\zcRefTypeSetup{definition}{cap}
\newtheorem{remark}[theorem]{Remark}
\AddToHook{env/remark/begin}{%
\zcsetup{countertype={theorem=remark}}}
\zcRefTypeSetup{remark}{cap}

\AddToHook{env/claim/begin}{%
\zcsetup{countertype={theorem=claim}}}
\zcRefTypeSetup{claim}{cap}

\zcRefTypeSetup{section}{cap}
\zcRefTypeSetup{subsection}{cap}
\zcRefTypeSetup{subsubsection}{cap}
\zcRefTypeSetup{figure}{cap}
\zcRefTypeSetup{item}{name-sg= ,name-pl= }
\zcRefTypeSetup{equation}{name-sg= ,name-pl= }

\newcommand\Gcb[3]{\Gamma(X,\varphi,#1,#2,#3)}
\newcommand\f\frac
\renewcommand\epsilon\varepsilon
\newcommand\intd{\mathop{}\!\mathrm{d}}
\newcommand\R{\mathbb R}

\newcommand\N{\mathbb N}
\renewcommand\phi\varphi
\newcommand\nFF[2]{\projX_{\#}\int_{\Gamma(X)\times[0,1]}#1 \intd(#2 \times \olm)}
\newcommand\diffmap[1]{#1'}
\newcommand\extd{\mathcal E}
\newcommand\phidiff{\diffmap{(\varphi\circ\gamma)}}
\newcommand\phidiffabs{|\phidiff|}
\newcommand\cext[1]{\overline{#1}}
\newcommand\idm{\operatorname{Id}}
\newcommand\II P
\DeclareMathOperator\Div{div}
\DeclareMathOperator{\spt}{spt}

\newcommand\ind[1]{1_{#1}}
\newcommand{\calH}{\mathcal{H}}
\newcommand\lm[1]{\mathcal L^n(#1)}

\newcommand{\olm}{\calH^1}
\newcommand{\TT}{\mathbf{T}}
\newcommand{\Ss}{\mathbb{S}}

\newcommand\mur{\nu}
\newcommand\oldnu{\rho}
\newcommand\restrict[2]{#1{\mathbin\upharpoonright}_{#2}}
\newcommand\ft[1]{\widehat{#1}}
\newcommand\ftl[1]{\mathcal F(#1)}

\newcommand{\dom}[1]{\operatorname{dom}#1}
\newcommand{\curveg}[1]{{#1}_{\varepsilon,R}}

\newcommand{\graph}{\operatorname{graph}}
\newcommand{\projX}{\pi}
\renewcommand{\P}{\eta}

\newcommand{\meas}{\mathcal{M}}

\title{Alberti representations and rectifiability of metric spaces}
\author{David Bate}
\author{Julian Weigt}
\affil{University of Warwick}
\affil{
\href{mailto:david.bate@warwick.ac.uk}{david.bate@warwick.ac.uk},
\href{mailto:julian.weigt@warwick.ac.uk}{julian.weigt@warwick.ac.uk}
}
\date{\today}

\begin{document}

\maketitle

\begin{abstract}
	We give a sufficient condition for a Borel subset $E\subset X$ of a complete metric space with $\calH^n(E)<\infty$ to be $n$-rectifiable.
	This condition involves a decomposition of $E$ into rectifiable curves known as an \emph{Alberti representation}.
	Precisely, we show that if $\restrict{\calH^n
	}{E}$ has $n$ independent Alberti representations, then $E$ is $n$-rectifiable.
	This is a sharp strengthening of prior results of Bate and Li.
	It has been known for some time that such a result answers many open questions concerning rectifiability in metric spaces, which we discuss.
\end{abstract}

\begingroup
\begin{NoHyper}%
\renewcommand\thefootnote{}\footnotetext{%
}%
\addtocounter{footnote}{-1}%
\end{NoHyper}%
\endgroup

\section{Introduction}

A Borel subset $E$ of a metric space $X$ is \emph{$n$-rectifiable} if there exist countably many $A_i\subset \R^n$ and Lipschitz $f_i\colon A_i\to X$ such that
\begin{equation*}
	\calH^n\biggl(E\setminus\bigcup_{i\in\N} f_i(A_i)\biggr)=0.
\end{equation*}
Rectifiable sets provide a fundamental notion of regularity that appears throughout analysis and differential geometry.
This notion has been extensively studied in the classical setting, that is when $X$ is some Euclidean space.
The properties of rectifiable sets, and sufficient conditions for rectifiability, are characterised in various structural theorems of geometric measure theory.
Equivalently, such structural theorems characterise \emph{purely $n$-unrectifiable} sets: Borel sets $S\subset X$ for which $\calH^n(S\cap E)=0$ for each $n$-rectifiable $E\subset X$.

Accompanying the recent interest in studying analysis and geometry in non-Euclidean metric spaces, there is a growing need for structural theorems for rectifiable sets in a metric space.
Kirchheim \cite{kirchheim} gives a generalisation of Rademacher's theorem for Lipschitz $f\colon A\subset \R^n\to X$ and, using this, generalises many of the results of classical geometric measure theory concerning rectifiable sets to metric spaces.
One consequence in particular is that the $f_i$ in the definition of a rectifiable set may be chosen to be bi-Lipschitz.

This article is a contribution to the study of sufficient conditions for rectifiability in a metric space, or equivalently, geometric properties of purely unrectifiable subsets of a metric space.

The conditions we consider are in terms of disintegrations of a measure into curves in the form of \emph{Alberti representations}.
The set $\Gamma(X)$ of curve fragments in $X$ consists of all 1-Lipschitz maps $\gamma$ defined on a compact subset $\dom\gamma\subset[0,1]$ into $X$.
After equipping this set with a suitable metric (see \zcref{curves}), an Alberti representation of a Borel measure $\mu$ on $X$ is a finite Borel measure $\P$ on $\Gamma(X)$ such that
\begin{equation}
	\label{AR-def-intro}
	\mu \ll \int_{\Gamma(X)}\restrict{\calH^1}\gamma \intd\P(\gamma)
	.
\end{equation}
We will also write \(\gamma\) for the image of \(\gamma\), so that \(\restrict{\calH^1}\gamma\) is the restriction of \(\calH^1\) to the image of \(\gamma\).
We say that several Alberti representations $\P_1,\ldots,\P_n$ are \emph{independent} if there exists a Lipschitz map $\varphi\colon X\to\R^n$ that distinguishes the curve fragments of each Alberti representation in the following way.
Up to a countable decomposition of $X$, there exist linearly independent (convex, one-sided) cones $C_1,\ldots,C_n \subset \R^n$ such that for each $1\leq i \leq n$ and $\P_i$-a.e.\ $\gamma\in\Gamma(X)$ we have
\begin{align*}
(\varphi\circ\gamma)'(t)
&\in
C_i\setminus \{0\}
&
\text{for $\calH^1$-a.e.\ }t&\in\dom\gamma
\end{align*}
(see \zcref{AR-independent}).

Fubini's theorem shows that Lebesgue measure on $\R^n$ possesses $n$ independent Alberti representations.
If $E\subset X$ is $n$-rectifiable, then by Kirchheim's bi-Lipschitz parametrisation theorem, we immediately see that $\restrict{\calH^n}{E}$ also possesses $n$ independent Alberti representations.
It is natural to ask if the converse to this statement holds, which we answer affirmatively.

\begin{theorem}
	\label{thm-alberti}
	Let $X$ be a complete metric space and let $E\subset X$ be Borel with $\calH^n(E)<\infty$ such that $\restrict{\calH^n}{E}$ has $n$ independent Alberti representations.
	Then $E$ is $n$-rectifiable.
\end{theorem}

This \zcref*[noref,nocap]{thm-alberti} is a strengthening of Bate and Li \cite[Theorem~1.2, (ii) $\Rightarrow$ (R)]{bateli}, where the same implication is shown, under the additional hypothesis that the lower $n$-dimensional density of $E$,
\begin{equation*}
	\Theta_*^n(E,x):= \liminf_{r\to 0} \frac{\calH^n(E\cap B(x,r))}{(2r)^n},
\end{equation*}
is positive at $\calH^n$-a.e.\ $x\in E$.

The results of \cite{bateli} have several important consequences in analysis on metric spaces.
It has been known for some time that proving \zcref{thm-alberti} in full generality answers several open questions concerning rectifiability in metric spaces.

\subsection{Consequences of \texorpdfstring{\zcref{thm-alberti}}{the main theorem}}\label{consequences}

The primary motivation to prove \zcref{thm-alberti} is its relationship to the following non-linear counter part to the classical Besicovitch--Federer projection theorem in any metric space.
\begin{theorem}
	\label{thm-BF}
	Let $X$ be a complete metric space and let $S\subset X$ be a purely $n$-unrectifiable set with $\calH^n(S)<\infty$.
	Then, for any $m\in\N$, a typical 1-Lipschitz map $f \colon X \to \R^m$  satisfies $\calH^n(f(S))=0$.

	Here ``typical'' means that the statement holds for a residual (in the Baire category sense) subset of all bounded 1-Lipschitz functions equipped with the supremum norm.
\end{theorem}
This theorem was first proven, under the additional assumption of positive lower density, in \cite[Theorem~1.1]{bateBF}.
The converse statement, that typical 1-Lipschitz images of rectifiable metric spaces have positive Hausdorff measure, is proven there in general without additional assumptions.
In some circumstances, a stronger converse holds \cite{jakub}.

To prove \zcref{thm-BF} one follows the strategy of \cite{bateBF}, but replacing \zcref{thm-density} with \zcref{thm-alberti}.
Specifically, \cite[Theorem~2.11]{bateBF} and \zcref{thm-alberti} show that any measure on a purely $n$-unrectifiable set with finite $\calH^n$ measure has an $(n-1)$-dimensional \emph{weak tangent field}.
\zcref{thm-BF} then follows from \cite[Theorem~5.4]{bateBF}.
A projection theorem for functions taking values in $\ell_\infty^m$ follows similarly, see \cite[Theorem~6.5]{bateBF}.

The non-linear projection theorem in \cite[Theorem~1.1]{bateBF} has consequences of its own.
By establishing \zcref{thm-BF} we are able to prove improved versions of these consequences.
Indeed, the argument of David and Le Donne \cite{ledonne} now gives the following.
\begin{theorem}
	\label{donnedavid}
	Let $X$ be a compact metric space with finite topological dimension $n$ and $\calH^n(X)<\infty$.
	Then $X$ contains an $n$-rectifiable subset of positive $\calH^n$ measure.
\end{theorem}
\zcref[S]{donnedavid} generalises Meier and Wenger \cite[Corollary~1.5]{wenger1} to all dimensions $n\geq 3$ and moreover removes the assumption of a geodesic metric.

Basso, Marti and Wenger \cite{wenger2} use the results of \cite{bateBF} in order to prove the rectifiability of certain metric manifolds.
Section 4 of \cite{wenger2} is required only to prove lower density bounds required for the non-linear projection theorem.

\zcref[S]{thm-alberti} establishes a rectifiability criterion of metric currents in the sense of Ambrosio and Kirchheim \cite{metric-currents}.
\begin{theorem}
	\label{currents}
	Let $\TT$ be a metric $n$-current in a complete metric space $X$ such that the support of $\|\TT\|$ has finite $\calH^n$ measure.
	Then $\TT$ is $n$-rectifiable.
\end{theorem}

This \zcref*[noref,nocap]{currents} strengthens \cite[Theorem~8.7]{metric-currents} from normal currents to all metric currents.
The classical case, for flat chains, is provided by White \cite[Corollary~6.1]{white}.
\zcref[S]{currents} follows from \zcref{thm-alberti} and the results of Schioppa \cite{schioppa-currents} that establish the existence of $n$ independent Alberti representations of the mass measure of a metric $n$-current.

The results in \cite{bateli} were partially motivated by characterising rectifiable subsets of \emph{Lipschitz differentiability spaces}, metric measure spaces that satisfy Cheeger's generalisation of Rademacher's theorem \cite{cheeger}.
\zcref[S]{thm-alberti} improves the results of \cite{bateli} to the following.
\begin{theorem}
	\label{thm-LDS}
	Any $n$ dimensional chart with finite $\calH^n$ measure in a Lipschitz differentiability space is $n$-rectifiable.
\end{theorem}
Indeed, this follows by combining \cite[Theorem~6.6]{structure} and \zcref{thm-alberti}.

Finally we remark that the original formulation of \cite{bateli} was in terms of a finite Borel measure $\mu$ with positive and finite lower and upper $n$-dimensional densities.
Along these lines, we see that \zcref{thm-alberti} holds for a general $\mu$ with positive upper $n$-dimensional density almost everywhere.
Indeed, such a measure is supported on a set of $\sigma$-finite $\calH^n$-measure.
The fact that $\mu\ll\calH^n$ follows from \cite{DPR}.

Since the first version of this article appeared, there have been several developments that make use of \zcref{thm-alberti}.
See \cite{2412.15794,2510.13525,2501.05974,2602.23322,2610.03968}.

\subsection{First steps in the proof of \texorpdfstring{\zcref{thm-alberti}}{the main theorem}}

In \cite[Theorem~2.18]{bateBF}, the first author had already observed a proof of \zcref{thm-alberti} when $X=\R^m$.
That argument used the work of De Philippis--Rindler \cite{DPR}, which establishes a structure theorem for measures satisfying a PDE constraint.
In particular, they show that any measure on $\R^n$ with $n$ independent Alberti representations is absolutely continuous with respect to Lebesgue measure.
Combined with the Besicovitch--Federer projection theorem, this implies that \zcref{thm-alberti} holds when $X=\R^m$.

That approach does not work in a general metric space since the non-linear projection theorem, \zcref{thm-BF}, is an aim we wish to \emph{conclude}.
Instead we prove \zcref{thm-alberti}, by showing show the hypotheses of \cite{bateli} are satisfied:
\begin{theorem}
	\label{thm-density}
	Let $X$ be a complete metric space and let $E\subset X$ be Borel with $\calH^n(E)<\infty$ such that $\restrict{\calH^n}{E}$ has $n$ independent Alberti representations.
	Then $\Theta_*^n(E,x)>0$ for $\calH^n$-a.e.\ $x\in E$.
\end{theorem}

From \cite{DPR} it is possible to also obtain a \emph{quantitative} result, slightly stronger than the qualitative structure theorem they state.
To prove \zcref{thm-density}, a first approach one my try is to use this refinement to deduce a \emph{quantitative} absolute continuity of measures on $\R^n$ with $n$ extremely regular Alberti representations, and then apply this result to the measure $\varphi_{\#}\restrict{\calH^n}{E\cap B(x,r)}$.
We precisely formulate \zcref{cor_pdeapplication} in this vein:
\begin{enumerate}
	\item \label{strong-hyp} 
The hypotheses, specifically the $L^1$ condition in \zcref{pdeapplication_l1bound}, is formulated as weakly as possible to avoid the immediate obstacles in applying it to $\varphi_{\#}\restrict{\calH^n}{E\cap B(x,r)}$ (see \zcref{need-L1}).
	\item The conclusion is useful to approach \zcref{thm-density}: for any $x$ and $r$ for which $\varphi_{\#}\restrict{\calH^n}{E\cap B(x,r)}$ satisfies the hypotheses,
	the conclusion implies $\calH^n(E\cap B(x,r)) \gtrsim r^n$.
\end{enumerate}
However, even taking into account \zcref{strong-hyp}, the hypotheses of \zcref{cor_pdeapplication} are too strong and it is not possible to show they are satisfied by $\varphi_{\#}\restrict{\calH^n}{E\cap B(x,r)}$ for $\calH^n$-a.e.\ $x\in E$ and all sufficiently small $r>0$.
Roughly, there are two fundamental issues that arise from the fact that $\restrict{\calH^n}{E}$ is not, a priori, asymptotically doubling.
Thus, a more sophisticated approach is required, in which \zcref{cor_pdeapplication} forms one part of the argument.

We will discuss our approach and the ideas behind it in \zcref{sec_approximate} once we have the necessary terminology to hand.
The article begins with two preliminary sections.
In \zcref{subsec_alberti} we recall the definition of an Alberti representation and related terminology.
In \zcref{sec_pde} we prove a quantitative absolute continuity result for measures on Euclidean space with finite divergence and deduce \zcref{cor_pdeapplication}.
In \zcref{three_lemmas} we will prove three lemmas that will be used in our proof.
The first gives a general construction of a measure defined as a disintegration into curve fragments, including a new way of defining a measure with an Alberti representation.
This method has the advantage of avoiding several measurability assumptions present in previous approaches and has a more general form.
The second provides a self improvement property for such measures.
In particular, it is possible to replace the absolute continuity in \zcref{AR-def-intro} by equality.

\subsection{Acknowledgements}

D.B.\ and J.W.\ were supported by the European Union's Horizon 2020 research and innovation programme (Grant agreement No.\ 948021),
in addition J.W.\ was funded by the European Union's Horizon research and innovation programme under the Marie Skłodowska-Curie Action 2023 No 101151034 (SRMF).
We would like to thank Tuomas Orponen for several discussions that led to \zcref{divergence}.

No generative AI was used in the preparation of this manuscript.

\section{Preliminaries}
\subsection{Alberti representations}
\label{subsec_alberti}
In this section we fix a metric space $X$.
We let $B(x,r)$ denote the closed ball centred on $x\in X$ of radius $r\geq 0$.
The set of all finite Borel regular (outer) measures on $X$ will be denoted $\meas(X)$.

\begin{definition}
    \label{curves}
    Define the set of \emph{curve fragments} $\Gamma(X)$ as the set of all 1-Lipschitz
    functions
	\begin{equation*}
		\gamma\colon \dom\gamma \subset [0,1] \to X
	\end{equation*}
    defined on a compact subset of $[0,1]$.
By \emph{curves} we understand curve fragments whose domain is an interval.

    For $\gamma\in\Gamma(X)$ let
    \[\graph(\gamma)=\{(t,\gamma(t)) \in [0,1]\times X: t\in\dom\gamma\}.\]
    We equip $\Gamma(X)$ with the Hausdorff distance between the graphs of curve fragments (treating the empty fragment as a discrete point).

The spaces $\Gamma(X)$ and $\Gamma(X) \times [0,1]$ are complete and separable and compact whenever $X$ is.
Moreover, the choice of metric on $\Gamma(X)$ implies that the map $(\gamma,t)\mapsto\gamma(t)$ is continuous on the closed subset of $\Gamma(X)\times[0,1]$ where it is defined.

For $\gamma\in\Gamma(X)$ and $t\in \dom\gamma$ we define its derivative by
\[
|\gamma'|(t) = \lim_{\dom\gamma\ni s\to t} \frac{d(\gamma(s),\gamma(t))}{|s-t|}
.
\]
The limit exists in a.e.\ $t\in\dom\gamma$, and a $1$-Lipschitz extension $\cext\gamma$ has the same derivative as $\gamma$ in a.e.\ $t\in\dom\gamma$.
This can be shown by isometrically embedding $\gamma$ into $\ell_\infty$ and applying \cite[Theorem~4.1.6]{ambtilli} to the the piecewise linear extension of $\gamma$ to $[0,1]$.
\end{definition}

\begin{definition}
	\label{generalAR}
	An \emph{Alberti representation} of a $\mu\in\meas(X)$ is a \(\P\in\meas(\Gamma(X))\) such that $\P$-a.e.\ $\gamma$ is injective and
	\begin{equation}
		\label{AR-def}
		\mu
		\ll
		\int_{\Gamma(X)}
		\restrict{\calH^1}\gamma
		\intd\P(\gamma)
		.
	\end{equation}
\end{definition}

The right hand side equals $\nFF\P\extd$ for $\extd$ the above zero extension of $\extd(\gamma,t) = |\gamma'|(t)$.

\begin{remark}
	If $X$ is complete and separable then,
	by \zcref{layer-cake} below,
	we obtain an equivalent definition of an Alberti representation if we replace the absolute continuity with equality
	(up to a countable decomposition of $X$).
\end{remark}

We use the following terminology to distinguish \emph{different} Alberti representations.

\begin{definition}
	\label{cones}
	For $w\in\Ss^{n-1}$ and $\theta>0$ let
	$C(w,\theta)$ be the cone consisting of those $v\in\R^n$ satisfying $\langle v,w\rangle \geq (1-\theta^2/2)|v|$.
	By a \emph{cone} we will always refer to a set of this form.
Cones $C_1,\ldots,C_n\subset \R^n$ are said to be \emph{independent} if any \(v_1,\ldots,v_n\in\R^n\) with $v_i\in C_i\setminus\{0\}$ are linearly independent.
\end{definition}

\begin{definition}
\label{AR-independent}
    Alberti representations $\P_1,\ldots,\P_n$ are \emph{independent} if there exists a
    \(1\)-Lipschitz map $\phi\colon X \to\R^n$ and a countable Borel decomposition
    \[X=\bigcup_{k\in\N} X_k\]
    and, for each $k\in\N$, there exist independent cones $C_1,\ldots,C_n\subset \R^n$ such that
    for each $1\leq i\leq n$, $\P_i$-a.e.\ $\gamma\in\Gamma(X)$ and $\olm$-a.e.\ $t\in \gamma^{-1}(X_k)$, we have $(\phi\circ\gamma)'(t)\in C_i\setminus \{0\}$.
\end{definition}

We also require quantitative control over aspects of independent Alberti representations.

\begin{definition}
For $\delta>0$, denote by $\Gamma(X,\delta)$ the set of all curve fragments $\gamma\in\Gamma(X)$ that are $(\delta,1)$-bi-Lipschitz, i.e.\ for any $t_1,t_2\in\dom\gamma$ we have \(d(\gamma(t_1),\gamma(t_2))\geq\delta|t_2-t_1|.\)

For a \(1\)-Lipschitz map $\varphi\colon X\to \R^n$, $e\in \Ss^{n-1}$, $\epsilon>0$ and $\delta>0$ we denote by $\Gcb e\varepsilon\delta$ the set of $\gamma\in \Gamma(X)$ with $\varphi\circ\gamma\in\Gamma(\varphi(X),\delta)$ and such that, for all $t_1,t_2\in\dom\gamma$ with \(t_1\leq t_2\),
\(
\varphi(\gamma(t_2))-\varphi(\gamma(t_1))
\in
C(e,\varepsilon)
.
\)

Note that $\Gamma(X,\delta)$ and $\Gcb{e}{\varepsilon}{\delta}$ are closed, and for $\gamma\in\Gcb{e_i}{\epsilon}{\delta}$ and almost every $t\in\dom\gamma$ we have
\begin{equation}
\label{phigammabilip}
\delta\leq|\phidiff(t)|\leq|\gamma'(t)|\leq1
.
\end{equation}
\end{definition}

We are able to \emph{refine} Alberti representations as follows.
The statement and proof of this result are a slightly more careful formulation of \cite[Corollary~5.9]{structure} and so we omit the proof.
\begin{lemma}
\label{filtering}
Suppose that \(\mu\in\meas(X)\) has $n$ independent Alberti representations with respect to a Lipschitz $\varphi\colon X\to\R^n$.
Then there exists a decomposition \(X=\cup_j X_j\) into Borel sets \(X_j\) such that for each \(j\in\N\) there are \(\delta_j,\tau_j>0\) for which the following holds:
For each \(\varepsilon>0\) there exists a finite Borel decomposition \(X_j=\cup_k X_j^k\) and, for each \(k\in\N\), a matrix \(\II=(e_1,\ldots,e_n)\) of unit vectors with \(|\II^{-1}|\leq\tau_j\) such that $\restrict{\mu}{X_j^k}$ has Alberti representation \(\P_1,\ldots,\P_n\) with \(\P_i\in\Gcb{e_i}{\varepsilon}{\delta_j}\).
\end{lemma}

\subsection{Regularity of measures \texorpdfstring{$L^1$-}{}close to invertible divergence-free measures}
\label{sec_pde}

In this \zcref*[noref,nocap]{sec_pde} we show a quantitative refinement of \cite[Theorem~1.1]{DPR} for the divergence operator and deduce \zcref{cor_pdeapplication} briefly mentioned in the introduction.
We use the essence of the proof from \cite{DPR}, but we refine each step in order to reach the quantitative conclusion.
It is perhaps possible to obtain a similar statement directly from \cite[Theorem~1.1]{DPR} using compactness arguments. We give a direct proof for brevity.

For $m\in\N$, $\meas(X;\R^m)$ will denote the set of $\R^m$ valued measures $\TT$ on $X$ with total variation $|\TT|\in \meas(X)$ and we write \(\|\TT\|\) for its total variation.
For $p>1$ let $p'=p/(p-1)$ denote the Hölder conjugate of \(p\).
By \(f\lesssim_ag\) we mean that for every \(a\) there exists a constant \(C\) such that \(f\leq Cg\).
If the inequality concerns Euclidean space then \(C\) may always depend on the dimension without specifying it in \(\lesssim_a\).

We require a preliminary lemma.

\begin{lemma}
\label{manyboundstol1}
Let \(1\leq p,q\leq\infty\) and \(b\in L^{1,\infty}(\R^n)\) with \(b^-\in L^p(\R^n)\) and \(\ft b\in L^q(\R^n)\).
Then
\(
\|b\|_{L^1(B(0,1))}
\lesssim_{p,q}
\|b\|_{1,\infty}^{\f1{p'}}
\|b^-\|_p^{\f1p}
+
\|\ft b\|_q
.
\)
\end{lemma}

\begin{proof}
	Take a nonnegative Schwartz function \(\psi\) with
	\[
	\ind{B(0,1)}
	\leq
	\psi
	\leq
	2\cdot\ind{B(0,2)}
	\]
	so that
	\[
	\|b\|_{L^1(B(0,1))}
	\leq
	\int
	\psi|b|
	=
	\int\psi b
	+
	2\int\psi b^-
	.
	\]
	Since \(\ft\psi\in L^{q'}(\R^n)\) we have
	\[
	\int\psi b
	=
	\int\ft\psi\ft b
	\leq
	\|\ft b\|_q
	\|\ft \psi\|_{q'}
	\lesssim_q
	\|\ft b\|_q
	\]
	and so it remains to bound \(\int\psi b^-\).
	We have
	\begin{align*}
	&
	\f12
	\int\psi b^-
	\leq
	\int_{B(0,2)}
	b^-
	\leq
	\int_0^{\min\{\|b\|_{1,\infty},\|b^-\|_p\}}
	\lm{B(0,2)}
	\intd\lambda
	\\
	&\qquad
	+
	\int_{\min\{\|b\|_{1,\infty},\|b^-\|_p\}}^{\|b^-\|_p}
	\lm{\{|b|>\lambda\}}
	\intd\lambda
	+
	\int_{\|b^-\|_p}^\infty
	\lm{\{|b^-|>\lambda\}}
	\intd\lambda
	.
	\end{align*}
	The first summand is equals \(\lm{B(0,2)}\) times
	\[
	\min\{\|b\|_{1,\infty},\|b^-\|_p\}
	\leq
	\|b\|_{1,\infty}^{\f1{p'}}\|b^-\|_p^{\f1p}
	.
	\]
	The second summand vanishes if \(\|b^-\|_p\leq\|b\|_{1,\infty}\), and otherwise equals
	\begin{align*}
	\int_{\|b\|_{1,\infty}}^{\|b^-\|_p}
	\lm{\{|b|>\lambda\}}
	\intd\lambda
	&\leq
	\int_{\|b\|_{1,\infty}}^{\|b^-\|_p}
	\frac{\|b\|_{1,\infty}}\lambda \intd\lambda
	\\
	&=
	\|b\|_{1,\infty}
	\log(\|b^-\|_p/\|b\|_{1,\infty})
	\lesssim_p
	\|b\|_{1,\infty}^{\f{1}{p'}}\|b^-\|_p^{\f1p}
	.
	\end{align*}
	Finally, we can bound the third summand by
	\begin{align*}
	\int_{\{|b|>\|b^-\|_p\}}
	b^-
	&\leq
	\lm{\{|b|> \|b^-\|_p\}}^{\f1{p'}}\|b^-\|_p
	\leq
	\|b\|_{1,\infty}^{\f1{p'}}\|b^-\|_p^{\f1p}
	.
	\end{align*}
\end{proof}

Consider \(\TT=(\TT_1,\ldots,\TT_n)\in \meas(\R^n;\R^{n\times n})\).
The divergence of $\TT$ is defined, row-wise, by duality with $\R^n$ valued Schwartz functions.

\begin{proposition}
\label{divergence}
Let \(\mur\in \meas(B(0,1))\) and let \(\TT\in\meas(\R^n;\R^{n\times n})\) be such that $\Div\TT \in \meas(\R^n;\R^n)$.
For any $1 < p < \f n{n-1}$,
\(\mur\) can be decomposed as
\[\mur=g+b\]
with \(g\in L^p(\R^n)\) and \(b\in\meas(\R^n)\) such that
\begin{align}
\label{divergence_good}
		\|g\|_p
&\lesssim_p
\|\mur\|+\|\Div\TT\|
		,\\
\label{divergence_bad}
\|b\|
&\lesssim_p
(\|\mur\|+\|\Div\TT\|)^{\f1p}
\|\idm\mur-\TT\|^{\f1{p'}}
,
\end{align}
where \(\idm\) denotes the identity matrix.
\end{proposition}

\begin{proof}
	By mollifying \(\mur\) and \(\TT\) and a standard use of Banach-Alaoglu  and Mazur, it suffices to consider the case that \(\mur,\TT\) are smooth functions.
	Write
	\[
		\Div(\idm\mur) = \Div(\idm\mur - \TT) + \Div\TT
	\]
	and take Fourier transforms
	\[
		\ft\mur(\xi)\xi = \ftl{\idm\mur - \TT}(\xi)\xi + \ft{\Div\TT}(\xi).
	\]
	Multiplying by \(\xi^\ast\) and adding $\ft\mur(\xi)$ leads to
	\[
		\ft\mur(\xi)(1 + \xi^2)
		=
		\xi^{\ast}\ftl{\idm\mur - \TT}(\xi)\xi + \ft\mur(\xi) + \xi^\ast\ft{\Div\TT}(\xi)
	\]
	and dividing by $1 + \xi^2$ gives
	\[
		\ft\mur(\xi)
		=
		\f{\xi^{\ast}\ftl{\idm\mur - \TT}(\xi)}{1+\xi^2}
		+
		\f{\ft\mur(\xi)}{1+\xi^2}
		+
		\f{\xi^\ast\ft{\Div\TT}(\xi)}{1+\xi^2}
		.
	\]

	Set \(b_0\) to be the inverse Fourier transform of the first summand and \(g_0\) the inverse Fourier transform of the second and third summand.
	Note that the Bessel potential \(B\) of order 1 belongs to $L^p$ and \(\ft B(\xi)=(1+\xi^2)^{-1/2}\) (\cite[6.1.2]{modernHA}).
	Thus, standard Calderon-Zygmund multiplier theory (\cite[Theorem~5.2.7]{classicalHA}) gives \zcref{divergence_good} for \(g_0\),
	\begin{align}
		\label{divergence_bweakbound}
		\|b_0\|_{1,\infty}
		&\lesssim
		\|\idm\mur-\TT\|
		,&
		\|\ft b_0\|_\infty
		&\lesssim
		\|\idm\mur-\TT\|
		.
	\end{align}
	By \zcref{manyboundstol1}, \zcref{divergence_bweakbound} and \zcref{divergence_good} for \(g_0\) we can conclude
	\begin{equation}
	\label{b0bound}
	\|b_0\|_{L^1(B(0,1))}
	\lesssim_p
	\|b_0^-\|_p^{\f1p}
	\|\idm\mur - \TT\|^{\f{1}{p'}}
	+
	\|\idm\mur - \TT\|
	.
	\end{equation}

	Recall that \(\mur\) is supported on \(B(0,1)\) and nonnegative.
	That means \(b_0^-\leq g_0\) and \(g_0^-\leq b_0\) and thus for
	\begin{align*}
	g&=(g_0^+-b_0^-)\ind{B(0,1)}
	,&
	b&=(b_0^+-g_0^-)\ind{B(0,1)}
	\end{align*}
	we have
	\begin{align*}
	0\leq g&\leq g_0
	,&
	0\leq b&\leq b_0
	,&
	\mur&=g+b
	.
	\end{align*}
	Since \(g_0\) satisfies \zcref{divergence_good}, so does \(g\).
	From \zcref{b0bound}, \zcref{divergence_good} for \(g_0\) and the fact that \(|b|\leq|b_0|\ind{B(0,1)}\) and \(|b_0^-|\leq|g_0|\), it follows that
	\[
	\|b\|_1
	\lesssim_p
	(\|\mur\|+\|\Div\TT\|)^{\f1p}
	\|\idm\mur - \TT\|^{\f{1}{p'}}
	+
	\|\idm\mur - \TT\|
	.
	\]
	If \(\|\idm\mur-\TT\|\leq\|\mur\|\) then this implies \zcref{divergence_bad}.
	If \(\|\idm\mur-\TT\|>\|\mur\|\) then \zcref{divergence_good,divergence_bad} hold for $(g,b)=(0,\mur)$.
\end{proof}

\begin{corollary}
\label{cor_pdeapplication}
	For any \(\tau,D>0\) there exist $c,d>0$ such that the following holds.
	Let \(r>0\), let \(\mur\in \meas(B(0,r))\) and let \(\TT\in\meas(\R^n,\R^{n\times n})\) be such that $\Div\TT \in \meas(\R^n,\R^n)$.
	Let \(\II\) be an invertible \(n\times n\)-matrix with \(|\II^{-1}|\leq\tau\) and assume
	\begin{equation}
	\label{pdeapplication_l1bound}
		\|\II\mur-\TT\|
		\leq
		d
		\|\mur\|
	\end{equation}
	and
	\begin{equation}
	\label{pdeapplication_divbound}
		r
		\|\Div\TT\|
		\leq
		D
		\|\mur\|
		.
	\end{equation}
	Then
	\[
		\lm{\spt\mur}
		\geq
		c
		r^n
		.
	\]
\end{corollary}

\begin{proof}
	Note that \(1/\det\II=\det(\II^{-1})\leq|\II^{-1}|^n\) for $|\II^{-1}|$ the matrix operator norm.
	Define \(\tilde\II=r|\II^{-1}|\II\)
and transform \(\tilde\mur=\mur\circ\tilde\II\) and \(\tilde\TT=\II^{-1}\TT\circ\tilde\II\).
Then \(|\tilde\II^{-1}|=1/r\) which ensures that \(\tilde\mur\) and \(\tilde\TT\) are supported on \(B(0,1)\).

Arbitrarily fix \(1<p<\f n{n-1}\).
By \zcref{divergence} we can decompose \(\tilde\mur=\tilde g+\tilde b\) with \(\tilde g,\tilde b\) satisfying \zcref{divergence_good,divergence_bad} in terms of \(\tilde\mur\) and \(\tilde\TT\).
That is, for \(g=\tilde g\circ\tilde\II^{-1}\) and \(b=\tilde b\circ\tilde\II^{-1}\), we have \(\mur=g+b\) which satisfy (after unpacking the scaling)
	\begin{align*}
		\|g\|_p
		&\lesssim_p
		r^{-n(1-1/p)}
		(
		\|\mur\|+r\tau\|\Div\TT\|
		)
		\leq
		(1+D\tau)
		r^{-n(1-1/p)}
		\|\mur\|,
		\\
		\|b\|
		&\lesssim_p
		((1+D\tau)\|\mur\|)^{\f1p}
		(d\tau\|\mur\|)^{\f1{p'}}
		=
		\tau(\tau^{-1}+D)^{\f1p}
		d^{\f{p-1}p}
		\|\mur\|
		.
	\end{align*}
	That means for some \(d^{\f{p-1}p}\gtrsim_p\tau^{-1}(\tau^{-1}+D)^{\f{-1}p}\) we have
	\(
		\|b\|
		\leq
		\|\mur\|/2
		.
	\)
	By Hölder's inequality we have
	\[
	\|\mur\|/2
	\leq
	\|g\|
	\leq
	\lm{\spt g}^{\f{p-1}p}
	\|g\|_p
	\leq
	\lm{\spt \mur}^{\f{p-1}p}
	\|g\|_p
	.
	\]
	Thus we find \(c\gtrsim_p(1+D\tau)^{\f{-p}{p-1}}\) for which the conclusion holds.
\end{proof}

\section{Three lemmas}
\label{three_lemmas}

We give a general definition of how to define a measure as an integral combination of one dimensional measures.

\begin{lemma}
\label{disintegrationmeasure}
Let $X$ be a metric space,
$\P\in \meas(\Gamma(X))$
and
\(\extd\colon\Gamma(X)\times[0,1]\to\R^n\)
a Borel function.
Let $\projX\colon \Gamma(X)\times [0,1]\to X$ be given by $\projX(\gamma,t)=\gamma(t)$.
The formula
\begin{equation*}
	\nFF\extd\P
	=
	\projX_{\#}[\extd\cdot(\P\times\olm)]
\end{equation*}
defines an element of \(\meas(X,\R^n)\).
Moreover, for any Borel $B\subset X$,
\begin{align}
	\label{rep-formula1}
	\nFF\extd\P(B)
	&=
	\int_{\Gamma(X)} \int_{\projX^{-1}(B)}
	\extd(\gamma,t)
	\intd\olm(t)\intd\P(\gamma)
	\\&=
	\notag
	\int_{\Gamma(X)}
	\gamma_{\#}[\extd(\gamma,\cdot) \olm](B)
	\intd\P(\gamma)
\end{align}
and
\begin{equation*}
	\label{measurability-of-rep}
	\gamma\mapsto
	\int_{\projX^{-1}(B)}
	\extd(\gamma,t)
	\intd\olm(t)
\end{equation*}
is Borel measurable.
\end{lemma}

\begin{proof}
Since $\R$ is separable, $\P\times\olm \in \meas(\Gamma(X)\times[0,1])$.
The subsequent statements follow from Fubini's theorem.
\end{proof}

We will use this construction for two choices of $\extd$.
For $\varphi\colon X\to \R^n$ Lipschitz, $\extd$ given by $\extd(\gamma,t) = (\varphi\circ\gamma)'(t)$ whenever the derivative exists and \(0\) otherwise, is Borel, see \cite[Lemma~2.8]{structure}.
Similarly, the zero extension of $\extd(\gamma,t) = |\gamma'|(t)$ is also Borel.
Further, if $\gamma\in\Gamma(X)$ is injective then
\[\int_{\gamma^{-1}(B)} |\gamma'|\intd\olm = \calH^1(\gamma \cap B),\]
see \cite[Theorem~4.4.2]{ambtilli}.
Thus, the measure appearing on the right hand side of \zcref{AR-def} is of the above form.

We will use the following self improvement property.

\begin{lemma}
	\label{layer-cake}
	Let $X$ be a complete and separable metric space, $\P\in\meas(\Gamma(X))$ and $\extd\colon \Gamma(X)\times [0,1]\to \R^n$ Borel.
	Suppose that $\mu\in\meas(X)$ is absolutely continuous with respect to the measure
	\[\nu=\nFF{\extd}{\P}.\]
	Suppose also that, for any $\gamma\in\Gamma(X)$, compact $K\subset \dom\gamma$  and almost every $t\in K$,  $\extd(\restrict{\gamma}{K},t)=\extd(\gamma,t)$.

	Then there exists a countable Borel decomposition $X=N\cup \bigcup_j X_j$ where $\mu(N)=0$ and each $X_j$ is compact
	such that
	\begin{equation*}
		\restrict\mu {X_j} = \nFF{\extd}{\P_j}
	\end{equation*}
	for each $j\in\N$.

Moreover, suppose \(F\subset\Gamma(X)\) is such that, for any $\gamma\in F$ and compact $K\subset \dom\gamma$, $\restrict{\gamma}{K}\in F$.
Then for each $j\in\N$, $\P(\Gamma(X)\setminus F)=0$ implies $\P_j(\Gamma(X)\setminus F)=0$.
\end{lemma}

\begin{proof}
	Since $\P$ is finite and $\extd$ is real valued, $\nu$ is $\sigma$-finite.
	After countably decomposing $X$, we may suppose both $\mu$ and $\nu$ are finite.
	Then by the Radon-Nikodym theorem, there exists a $g\in L^1(\nu)$ such that $\mu= g\nu$.
	By Lusin's theorem, there exist countably many compact sets $X_j\subset X$ on which $g$ is continuous such that $\mu(X\setminus \cup_j X_j)=0$.
	Fix $j\in \N$ and let $m_j=\max_{X_j}g$.

	The map $\kappa_j:\gamma\mapsto\gamma\cap X_j$ is Borel (\cite[Theorem 18.10]{zbMATH05265624}).
	Since $g$ is continuous, $\lambda\mapsto g^{-1}([\lambda,m_j])\cap X_j$ is continuous.
	Thus, the map $p_j:\Gamma(X)\times[0,m_j]\rightarrow\Gamma(X_j)$ given by $p_j(\gamma,\lambda) =\gamma \cap g^{-1}([\lambda,m_j])\cap X_j$ is Borel.
	Define $\P_j= {p_j}_{\#} (\P\times \restrict{\olm}{[0,m_j]})$, so that the moreover statement is satisfied.

Finally, let $B\subset X_j$ be Borel.
	Then by \zcref{rep-formula1}, Fubini and
	the the assumption on $\extd$ in the second equality,
	\begin{align*}
		\mu(B) &=
		\int_0^{m_j}
		\nu(\{x\in B : g(x) \geq \lambda\})
		\intd\lambda
		\\&=
		\int_0^1
		\int_{\Gamma(X)}
		\int_0^{m_j}
		\ind{\{\gamma(t)\in B,\ g(\gamma(t)) \geq \lambda\}}
		\extd(\gamma,t)
		\intd\P(\gamma)\intd\lambda\intd t
		\\&=
		\int_0^1\int_{\Gamma(X)}\int_0^{m_j}
		\ind{\{p_j(\gamma,\lambda)(t)\in B\}}
		\extd(p_j(\gamma,\lambda),t)
		\intd \P(\gamma)\intd\lambda\intd t
		\\&=
		\int_0^1\int_{\Gamma(X)}
		\ind{\{\gamma(t)\in B\}}
		\extd(\gamma,t)
		\intd \P_j \intd t
		\\&=
		\biggl[
		\nFF{\extd}{\P_j}
		\biggr](B)
		.
	\end{align*}
\end{proof}

The moreover statement of the preceding result implies that the properties of an Alberti representation introduced in \zcref{subsec_alberti} (direction, independence, etc.) are preserved under the self improvement property.

We need to localise a measure in a metric space to small balls.
To do this whilst controlling the divergence of the resulting measure, we use a Lipschitz cutoff function and the following lemma.
\begin{lemma}
	\label{div-psi}
	Let $X$ be a metric space.
	For $x\in X$ and $r>0$ let $\psi\colon X\to \R$ be $\frac{1}{r}$-Lipschitz with $\spt\psi\subset B(x,2r)$.
	Let $\gamma\in\Gamma(X,\delta)$ with $\dom\gamma=[a,b]$ and $\gamma(a),\gamma(b)\not\in B(x,2r)$.
	Then
	\begin{equation}
		\label{div-calc}
		\|\Div(\varphi_\#(\psi\circ\gamma\cdot(\varphi\circ\gamma)'\olm))\|
	\leq
	\frac{1}{\delta r}\int_{\gamma^{-1}(B(x,2r))}|(\varphi\circ\gamma)'(t)|\intd t
	.\end{equation}
\end{lemma}

\begin{proof}
	Let $u\colon \R^n\to\R$ smooth with $\|u\|_\infty\leq 1$.
	The chain rule and integration by parts give
\begin{align*}
	\int
	\nabla u
	\intd(\varphi_\#(\psi\circ\gamma\cdot(\varphi\circ\gamma)'\olm))
	&=
	-
	\int
	(\nabla u)(\varphi(\gamma(t)))\psi\circ\gamma(t)(\varphi\circ\gamma)'(t)\intd t
	\\
    	&=
	-
	\int(u\circ(\varphi\circ\gamma))'(t)\psi\circ\gamma(t)\intd t
	\\
	&=
	\int(u\circ(\varphi\circ\gamma))(t) (\psi\circ\gamma)'(t)\intd t
	.
\end{align*}
Taking the supremum over such \(u\) yields
\begin{align*}
\|\Div(\varphi_\#(\psi\circ\gamma\cdot(\varphi\circ\gamma)'\olm))\|
&\leq
\int|(\psi\circ\gamma)'(t)|\intd t
\\
&\leq
\frac{1}{r}\int_{\gamma^{-1}(B(x,2r))}|\gamma'(t)|\intd t
\\
&\leq
\frac{1}{\delta r}\int_{\gamma^{-1}(B(x,2r))}|(\varphi\circ\gamma)'(t)|\intd t.
\end{align*}
\end{proof}
\section{The main argument}
\label{sec_approximate}

We begin by discussing the initial idea of this section.
Although the discussion is longer than the formal proof, there are several subtle points that create difficulty later on, which are worth emphasising.
Until \zcref{curvewiseRnbound}, we will use the following notation to aid discussion.

Let $X$ be a normed vector space, $\P\in\meas(\Gamma(X))$ and
\[
	\mu=\int_{\Gamma(X)}\restrict{\calH^1}\gamma \intd\P(\gamma)
	.
\]
For $x\in X$ and $r>0$, in order to apply \zcref{cor_pdeapplication}, we wish to construct a measure $\P'\in\meas(\Gamma(X))$ supported on \emph{curves} -- curve fragments whose domain is an interval -- such that
\[
\nu
=
\int_{\Gamma(X)}\restrict{\calH^1}\gamma \intd\P'(\gamma)
\]
satisfies
\begin{equation}
	\label{mu-nu-diff}
	|\mu-\nu|(B(x,r)) < \epsilon \mu(B(x,r)).
\end{equation}
It is natural to attempt to do this by taking the piecewise linear extension of each $\gamma$.

\begin{definition}
	\label{def_curveext}
	For a normed vector space space $X$ and $\gamma\in\Gamma(X)$ define $\cext{\gamma}$ to be the piecewise linear extension of $\gamma$ to the convex hull of $\dom\gamma$.

It is readily seen that the map $\gamma\mapsto \cext{\gamma}$ is Lipschitz and maps $\Gamma(X,\delta)$ and $\Gcb e\varepsilon\delta$ into themselves.
In particular, if $\gamma\in\Gamma(X,\delta)$ and \(\gamma(t)\in B(x,2r)\) then
\begin{equation}
\label{deltaforcescurvetoleave}
\cext\gamma^{-1}(B(x,2r))
\subset
[t-4r/\delta,t+4r/\delta]
.
\end{equation}
\end{definition}

If one takes the linear extension of \emph{all} fragments intersecting $B(x,r)$, it is impossible to bound the difference of $\mu$ and $\nu$.
Instead we will only extend curve fragments with high Lebesgue density in $B(x,r)$.

\begin{definition}\label{goodpoints-large}
	Let $X$ be a metric space.
For $\gamma\in\Gamma(X)$, $R>0$ and $0<\varepsilon<1$, let $\curveg\gamma$ be $\gamma$ restricted to those $t\in\dom\gamma$ for which $[t-R,t+R]\subset[0,1]$ and
\[\olm([t-r,t+r]\cap\dom\gamma)\geq 2(1-\epsilon)r\]
for all $0\leq r\leq R$.

Since $K\mapsto \calH^1(K)$ is upper semicontinuous with respect to Hausdorff convergence on $\R$, the set $\{(\gamma,t):\gamma\in\Gamma(X),\ t\in\dom{\curveg\gamma}\}$ is closed.
\end{definition}

For the moment assume $\P\in\meas(\Gamma(X,1/2))$.
For $\epsilon,R>0$ and $\gamma\in\Gamma(X,1/2)$, if $t\in\dom\curveg\gamma$ with $\gamma(t)\in B(x,r)$ then
\[
	\olm(\cext{\gamma}\setminus \gamma \cap B(x,r)) \lesssim \epsilon r.
\]
By extending only curve fragments of this form when defining $\nu$, we obtain
\begin{equation}
\label{crude-bound}
|\mu-\nu|(B(x,r))\lesssim \epsilon r \P(\{\gamma\in\Gamma(X,1/2): \exists t\in\dom\curveg\gamma, \gamma(t)\in B(x,r)\}).
\end{equation}
Recall that we require a bound similar to \zcref{mu-nu-diff}.
To achieve this we use the Alberti representation of $\mu$ and a bound of the form $\olm(\gamma\cap B(x,r))\gtrsim r$ for the fragments appearing in the set on the right hand side of \zcref{crude-bound}.
In general such a bound is not available,
for example if $\gamma$ intersects $B(x,r)$ only near the boundary of $B(x,r)$.
However, the following does hold
\begin{equation}
	\label{inter-lower}
	\olm(\gamma\cap B(x,2r))\gtrsim r.
\end{equation}
Combining \zcref{crude-bound,inter-lower} gives
\begin{equation*}
|\mu-\nu|(B(x,r))\lesssim \epsilon \mu(B(x,2r))
\end{equation*}
Notice that $\mu(B(x,2r))$ can be significantly larger than $\mu(B(x,r))$, which is required for \zcref{mu-nu-diff}.
This will cause difficulty later on.

In the preceding discussion we left the curve fragments $\gamma\in\Gamma(X,1/2)$ with $\dom\curveg\gamma \cap B(x,r) = \emptyset$ untouched.
In our true definition of $\nu$ we will remove those fragments, and use the following measure to control the resulting additional error:
\begin{definition}
Given $\epsilon,R>0$ and $\P\in\meas(\Gamma(X))$ define
\[
\oldnu_R
=
\projX_{\#}[\restrict{(\P\times\olm)}{\{
(\gamma,t):t\in\dom\gamma\setminus\dom{\curveg\gamma}
\}}]
.
\]
Although $\oldnu_R$ also depends on $\epsilon$ and $\eta$, this will be clear from the context and so we omit the dependence from the notation for brevity.
\end{definition}

The total measure of $\oldnu_R$ can be chosen arbitrarily small.
\begin{lemma}
	\label{total-small}
	For any metric space $X$, $\P\in\meas(\Gamma(X))$ and $\varepsilon>0$,
	$
    \|\oldnu_R\|
    \to 0
    $ as $R\to 0$.
	Moreover, if $\P\in\meas(\Gamma(X,\delta))$ for some $\delta>0$ then $\oldnu_R\leq \mu/\delta$.
\end{lemma}

\begin{proof}
	For fixed $\varepsilon>0$, the Lebesgue density theorem implies that, as $R\to 0$, $\olm(\dom\gamma\setminus\dom{\curveg\gamma}) \to 0$.
Thus,
	\[\{(\gamma,t):t\in\dom\gamma\setminus\dom{\curveg\gamma}\}\]
	monotonically increases to a set of full $\P\times\olm$ measure in \[\{(\gamma,t)\in\Gamma(X)\times[0,1]:t\in\dom\gamma\}.\]
That is,
	\[
    \lim_{R\rightarrow0}
		\|\rho_R\|
    =
    \lim_{R\rightarrow0}
    (\P\times\olm)(\{
(\gamma,t):t\in\dom\gamma\setminus\dom{\curveg\gamma}
\})
	=0
	\]
	The fact that $\oldnu_R\leq \mu/\delta$ follows from the fact that, in this case, $|\dot\gamma|(\delta)\geq \delta$ for $\P$-a.e.\ $\gamma$ and $\olm$-a.e.\ $t$.
\end{proof}

We now precisely formulate the ideas of the previous discussion.
Note, in the following lemma we consider the balls $B(x,2r)$ and $B(x,3r)$ rather than $B(x,r)$ and $B(x,2r)$, to account for a Lipschitz cutoff of $B(x,r)$ that will appear later.

\begin{lemma}
\label{curvewiseRnbound}
Let $X$ be a normed vector space space, $\delta>0$ and \(\P\in \meas(\Gamma(X,\delta))\).
For $R,\varepsilon>0$, $0<r<\delta R/4$, and $x\in X$ define $\mu,\nu,\in\meas(X)$ by
\begin{align*}
\mu
&=
\int_{\Gamma(X,\delta)} \restrict{\olm}{\gamma}\intd\P(\gamma)
,&
\nu
&=
\int_{\{
\gamma:
\curveg\gamma\cap B(x,2r)\neq\emptyset
\}}
\restrict{\olm}{\cext{\gamma}}\intd\P(\gamma)
.
\end{align*}
Then
\begin{align}
\label{smalldifference}
|\mu-\nu|(B(x,2r))
&\leq
\f{8\varepsilon}{\delta^2}
\mu(B(x,3r))
+
\rho(B(x,2r))
,\\
\label{notmuchmassadded}
\nu(B(x,2r)
&\leq
\Bigl(1+\frac{8\varepsilon}{\delta^2}\Bigr)
\mu(B(x,3r))
.
\end{align}
\end{lemma}

\begin{proof}
We divide $|\mu-\nu|(B(x,2r))$ into two terms,
\begin{align}
	\label{bound-diff-Gammag}
	|\mu-\nu|(B(x,2r))
	&\leq
	\int_{
\{
\gamma\in\Gamma(X,\delta):
\curveg\gamma\cap B(x,2r)\neq\emptyset
\}
}\olm(\cext{\gamma}\cap B(x,2r)\setminus \gamma)\intd\P(\gamma)
	\\
	&\quad+
	\int_{
\{
\gamma\in\Gamma(X,\delta):
\curveg\gamma\cap B(x,2r)=\emptyset
\}
}\olm(\gamma\cap B(x,2r))\intd\P(\gamma)
	\notag
\end{align}
Since $\gamma$ is 1-Lipschitz the second summand in \zcref{smalldifference} bounds the second summand in \zcref{bound-diff-Gammag},
and it remains to do the same for the respective first summands.

Let $\gamma\in\Gamma(X,\delta)$ and $t\in\dom{\curveg\gamma}$ with $\gamma(t)\in B(x,2r)$.
Then
\begin{equation}
\label{curvestaysinlargerball}
\olm(
\gamma \cap
B(x,3r)
)
\geq\delta
\olm(
[t-r,t+r]
\cap
\dom\gamma
)
\geq
2(1-\varepsilon)\delta r
\geq \delta r
.
\end{equation}
Further, since \(\gamma\) is \(1\)-Lipschitz and $4r/\delta<R$, \zcref{deltaforcescurvetoleave} implies
\begin{equation}
\label{extendedcurvewithoutcurve}
\olm(\cext{\gamma}\setminus \gamma \cap B(x,2r))
\leq
\olm([t-4r/\delta,t+4r/\delta]\setminus \dom\gamma)
\leq
8\epsilon r/\delta.
\end{equation}
Integrating \zcref{extendedcurvewithoutcurve,curvestaysinlargerball} over $\gamma$ proves \zcref{notmuchmassadded}, as well as bounds the first summand in \zcref{bound-diff-Gammag} by the first summand in \zcref{smalldifference}, finishing the proof of \zcref{smalldifference}.
\end{proof}

We now turn our attention to the upper bounds in \zcref{notmuchmassadded,smalldifference} which are too large: in order to apply \zcref{cor_pdeapplication}, we require upper bounds comparable to $\epsilon \mu(B(x,2r))$ and $\mu(B(x,2r))$ respectively.
We will \emph{not} show that such bounds hold for $\calH^n$-a.e.\ $x\in E$ all sufficiently small $r>0$.
In the following result we establish sufficient conditions for when it is possible to apply \zcref{cor_pdeapplication} and, under those conditions, show how to arrive to our desired conclusion.
Following this result, and for the remainder of the article, we focus on how to handle the case when these sufficient conditions are not satisfied.

\begin{proposition}
\label{pro_inductionstep}
For any \(\tau,\delta,C>0\) there exist \(\varepsilon,c>0\) such that the following holds:

Let $X$ be a metric space and let \(\varphi\colon X\to\R^n\) be 1-Lipschitz.
Let \(e_1,\ldots,e_n\) be unit vectors with \[|(e_1,\ldots,e_n)^{-1}|\leq\tau.\]
For each $1\leq i \leq n$ let $\P_i\in\meas(\Gcb{e_i}{\epsilon}{\delta})$ such that each of the measures
\[
\nFF{\phidiffabs}{\P_i}
\]
agree and denote the common measure by $\mu$.
For each \(1\leq i\leq n\) and $R>0$ define
\begin{align*}
\oldnu_{i,R}
&=
\projX_{\#}[\restrict{(\P_i\times\olm)}{\{
(\gamma,t):t\in\dom\gamma\setminus\dom{\curveg\gamma}
\}}]
.
\end{align*}

Then for every \(x\in X\) and \(0<r<\delta R/4\) with
\begin{equation}
\label{doublingassumption}
\mu(B(x,3r))
\leq
C
\mu(B(x,r))
\end{equation}
and
\begin{equation}
\label{eq_mostpointsgood}
\oldnu_{i,R}(B(x,2r))
<
C\varepsilon
\mu(B(x,r))
\end{equation}
we have
\[
\calH^n(
\spt(\varphi_\#(\restrict\mu{B(x,2r)}))
)
\geq
cr^n
.
\]
\end{proposition}

\begin{proof}
	Let $\iota\colon X\to V$ be a bi-Lipschitz embedding of $X$ into a Banach space.
	Then $\iota'(x):= (\varphi(x),\iota(x))$ is a bi-Lipschitz embedding into a Banach space such that the projection onto the first component agrees with $\varphi$ on the image of $X$.
	Since the conclusion is bi-Lipschitz invariant, we may suppose that $X$ is a Banach space and that $\varphi$ is linear.

	In this case, $\cext{\gamma} \in \Gcb{e_i}{\varepsilon}{\delta}$ for $\P_i$-a.e\ $\gamma$.
	Define $\cext{\P_i}$ to be the pushforward of
	$\restrict{\P_i}{\{
\gamma:
\curveg\gamma\cap B(x,2r)\neq\emptyset
\}}$ under the map $\gamma\mapsto \cext{\gamma}$, so that \[\cext{\P_i} \in \meas(\Gcb{e_i}{\varepsilon}{\delta}).\]
	Denote
	\[
	\nu_i=\nFF{\phidiffabs}{\cext{\P_i}}.
	\]
	Recall from \zcref{phigammabilip} that for $\gamma\in\Gcb{e_i}\varepsilon\delta$ and almost every $t\in\dom\gamma$ we have
	\(
	|\phidiff(t)|\leq1\leq|\phidiff(t)|/\delta
	.
	\)
	Thus, by \zcref[S]{curvewiseRnbound} we can conclude
	\begin{align}
		\notag
		|\mu- \nu_i|(B(x,2r))
		&\leq
		\biggl|\nFF{}{(\P_i-\cext{\P_i})}\biggr|(B(x,2r))
		\\
		\notag
		&\leq
		\frac{8\epsilon}{\delta^3} \mu(B(x,3r)) + \oldnu_{i,R}(B(x,2r))
		\\
		\label{pde-hyp1}
		&\leq
		(8/\delta^3+C)\epsilon \mu(B(x,3r))
		,\\
		\label{numu}
		\nu_i(B(x,2r))
		&\leq
		\Bigl(\frac1\delta+\frac{8\varepsilon}{\delta^3}\Bigr)
		\mu(B(x,3r))
		.
	\end{align}
	for each $1\leq i \leq n$, using \zcref{eq_mostpointsgood} for the final inequality.

	Define the Lipschitz cut off \(\psi\colon X\to [0,2]\) by
\[
    \psi(y)
    =
    \max\Bigl\{2-\frac{|y-x|}{r},0\Bigr\}
.
\]
Note that $\psi$ is $\frac{1}{r}$-Lipschitz and $\ind{B(x,r)}\leq\psi\leq 2\cdot \ind{B(x,2r)}$.
For each $1\leq i \leq n$ let
	\[T_i = \varphi_\#\biggl(\psi\nFF{\phidiff}{\cext{\P_i}}\biggr),\]
	a vector valued measure on $\R^n$.
	Since for $\cext{\P_i}\times\olm$-a.e.\ $(\gamma,t)$ we have
	\[
	\biggl|
	\frac{(\varphi\circ\gamma)'(t)}{|(\varphi\circ\gamma)'(t)|}
	-
	e_i
	\biggl|
	\leq
	\varepsilon
	,
	\]
	integrating yields
	\begin{equation}
	\label{smallangleconclusion}
	\|
	T_i - \varphi_\#(\psi\nu_i)e_i
	\|
	\leq
	\varepsilon\|\varphi_\#(\psi\nu_i)\|
	,
	\end{equation}

For any $\gamma\in\Gcb{e_i}{\epsilon}{\delta}$ with $\curveg\gamma\cap B(x,2r)\neq\emptyset$,
since $4\delta r<R$, \zcref{deltaforcescurvetoleave} implies that \zcref{div-calc} holds for $\cext\gamma$.
For each $1\leq i \leq n$, integrating with respect to $\cext{\P_i}$ gives
\begin{equation}
\|\Div(T_i)\|
\leq
\frac1{\delta r}\nu_i(B(x,2r))
.
\label{pde-hyp2}
\end{equation}

Taking $\varepsilon$ sufficiently small such that $\varepsilon/\delta$ is sufficiently small, we can conclude from \zcref{pde-hyp1,pde-hyp2,numu,smallangleconclusion,doublingassumption}
\begin{align*}
\|T_i-\varphi_\#(\psi\mu)e_i\|
&\leq
2C(9/\delta^3+C)\varepsilon
\|\varphi_\#(\psi\mu)\|
,\\
\|\Div(T_i)\|
&\leq
\frac{2C}{\delta^2 r}
\|\varphi_\#(\psi\mu)\|
,
\end{align*}
allowing us to apply \zcref{cor_pdeapplication}, which provides the required $c,\epsilon>0$.
\end{proof}

\begin{remark}
\label{need-L1}
In the preceding proof, the extended measures $\nu_i$ will typically not be mutually absolutely continuous.
This is precisely the reason why our formulation of \zcref{cor_pdeapplication} only requires the $L^1$ bound \zcref{pdeapplication_l1bound}.
There exist other quantitative results related to \cite{DPR}.
For example, in \cite{quant-PDE} a higher integrability result is obtained under the assumption of an $L^\infty$ hypothesis.
\end{remark}

The remainder of the article is devoted to understanding when \zcref{pro_inductionstep} may be applied (when \zcref{doublingassumption,eq_mostpointsgood} hold), and obtaining an alternate solution when it cannot.

To obtain \zcref{eq_mostpointsgood} for $\mu$-a.e.\ $x$ and all sufficiently small $r>0$, it would be necessary to use the \emph{global} property of small total measure given by \zcref{total-small} to obtain a \emph{local} property for most $x$.
Such an argument depends on the Lebesgue density theorem, which is not available to us.
Instead, we must find a way to proceed that only uses the following weaker conclusion.

\begin{lemma}
\label{lem_highdensityapprox}
Let \(E\subset X\) with \(\calH^n(E)<\infty\) and suppose that for each \(m\in\N\), \(\oldnu_m\in\meas(X)\) with \(\restrict{\calH^n}E\geq\oldnu_1\geq\oldnu_2\geq\ldots\geq0\) and \(\|\oldnu_m\|\rightarrow0\).
Then for \(\calH^n\)-almost every \(x\in X\) we have
\begin{equation}
\label{eq_highdensityapprox}
    \lim_{m\rightarrow\infty}\lim_{r\rightarrow0}
    \f{
    \oldnu_m(B(x,r))
    }{
    r^n
    }
    =0
    .
\end{equation}
\end{lemma}

\begin{proof}
By \cite[Theorem~2.10.18]{federer},
for \(\calH^n\)-almost every \(x\in X\setminus E\) we have \(\lim_{r\rightarrow0}\calH^n(E\cap B(x,r))/r^n=0\), so it remains to consider \(x\in E\).

Let \(\sigma_m\) be the Radon-Nikodym derivative of $\oldnu_m$ with respect to $\restrict{\calH^n}E$, so that the $\sigma_m$ monotonically decrease to 0 and are bounded above by 1.
For each \(k\in\mathbb{N}\) let \(m_k\in\N\) be such that
\[
N_k
:=
\{
x\in E:
\sigma_{m_k}(x)>2^{-k}
\}
\]
satisfies $\calH^n(N_k)<2^{-k}$.
For \(\calH^n\)-almost every \(x\not\in N_k\) there is an \(r_{x,k}>0\) such that, for all \(m\geq m_k\) and \(r<r_{x,k}\),
\[
\oldnu_m(B(x,r)\cap N_k)
\leq
\calH^n(B(x,r)\cap N_k)
\leq
2^{-k} r^n
.
\]
Further, by \cite[Theorem~2.10.17]{federer}, for \(\calH^n\)-almost every \(x\not\in N_k\) and \(m>m_k\), by reducing $r_{x,k}$ if necessary, for \(r<r_{x,k}\) we have
\[
\oldnu_m(B(x,r)\setminus N_k) \leq 2^{-k}\calH^n(B(x,r)\cap E) \leq 2^{-k}(3r)^n.
\]
Therefore,
\[
\oldnu_m(B(x,r))
=
\oldnu_m(B(x,r)\cap N_k)
+
\oldnu_m(B(x,r)\setminus N_k)
\leq
2^{-k}
(1+3^n)
r^n
.
\]
Define \(N=\bigcap_{k\in\mathbb{N}}\bigcup_{m\geq k}N_m\).
Then \(\calH^n(N)=0\) and for \(\calH^n\)-almost every \(x\in E\setminus N\) we have
\zcref{eq_highdensityapprox}.
\end{proof}

The conclusion of the previous \zcref*[noref,nocap]{lem_highdensityapprox} is much weaker than what is required to apply \zcref{pro_inductionstep}.
Indeed, deducing \zcref{eq_mostpointsgood} would first require us to show $\mu(B(x,r))\gtrsim r^n$ for almost every $x$ and all sufficiently small $r$, which is exactly the conclusion we seek!
In our approach, we instead make use of the properties of Hausdorff measure which provide a sequence of scales $r_i\to 0$ for which $\mu(B(x,r_i)) \gtrsim r_i^n$.
Together with a careful setup, this allows us to initialise a bootstrap argument based on the following dichotomy.
We may assume that a lower density bound holds at some scale, so that \zcref{eq_highdensityapprox} implies \zcref{eq_mostpointsgood}.
The \emph{presence} of the doubling condition \zcref{doublingassumption} allows us to apply \zcref{pro_inductionstep} and obtain a lower density bound at the next highest scale;
on the other hand, the \emph{absence} of \zcref{doublingassumption} directly allows us to improve the assumed lower density bound from the current scale to the next larger scale.

\begin{proposition}
\label{almosteverywheredensity}
For \(C=3^n\) and \(\tau,\delta>0\) obtain \(\varepsilon,c>0\) from \zcref{pro_inductionstep}.
Let $X,\varphi,e_i,\P_i$ and $\mu,\oldnu_{i,R}$ be as in \zcref{pro_inductionstep} and in addition assume that $\mu$ is of the form \(\mu=\restrict{\calH^n}E\) for $E\subset X$.
Then for \(\calH^n\)-almost every \(x\in E\) we have
\[
\Theta_{*}^{n}(E,x)
:=
\liminf_{r\rightarrow0}
\f{
\calH^n(E\cap B(x,r))
}{
(2r)^n
}
\geq
9^{-n}
\min\{c,2^n\}
.
\]
\end{proposition}

\begin{proof}
If \zcref{pro_inductionstep} yields a number $c>2^n$ then we set $c=2^n$ for the rest of the proof.
The conclusion of \zcref{pro_inductionstep} continues to hold for that value.

By \cite[Theorem~2.10.17]{federer},
\[
\limsup_{r\to 0} \frac{\mu(B(x,r))}{(2r)^n}
>
3^{-n}
\geq
6^{-n}c
\]
for $\mu$-a.e.\ $x$. That is, there is a sequence \(r_k\to0\) such that, for each \(k\in\N\),
\begin{equation}
\label{lowerdensityradius}
\mu(B(x,r_k))
\geq
c(r_k/3)^n.
\end{equation}

Recall from \zcref{phigammabilip} that
\(
|\phidiff(t)|\geq\delta
,
\)
which implies $\oldnu_{i,R}\leq\mu/\delta$.
Thus, \zcref{total-small,lem_highdensityapprox} imply that, for \(\mu\)-almost every \(x\in X\), there exists an \(R>0\) such that for all \(0<r<R\),
\begin{equation}
\label{eq_mostpointsgoodr}
\oldnu_{i,R}(B(x,2r))
\leq c\varepsilon
r^n
.
\end{equation}
We now fix an $x\in X$ for which \zcref{lowerdensityradius,eq_mostpointsgoodr} hold and the corresponding $r_k,R>0$.

For a moment fix $k$ with $r_{k}<\delta R/4$.
Note that \zcref{eq_mostpointsgoodr,lowerdensityradius} imply \zcref{eq_mostpointsgood}.
Thus we ask if \zcref{doublingassumption} holds.
That is, if
\begin{equation*}
\mu(B(x,3r_k))
\leq
3^n
\mu(B(x,r_k))
.
\end{equation*}
If so, we invoke \zcref{pro_inductionstep}.
Since \(\mu=\restrict{\calH^n}{E}\) and $\varphi$ is 1-Lipschitz,
we deduce
\[
\mu(B(x,2r_k))
\geq
\calH^n(\varphi(E\cap B(x,2r_k)))
=
\calH^n(\spt(\varphi_\#(\restrict\mu{B(x,2r_k)})))
\geq
cr_k^n
.
\]
On the other hand,
if \zcref{doublingassumption} fails, then \zcref{lowerdensityradius} implies
\begin{equation*}
\mu(B(x,3r_k))
\geq
3^n
\mu(B(x,r_k))
\geq 3^n c(r_k/3)^n = c r_k^n.
\end{equation*}
In either case,
\begin{equation*}
\mu(B(x,3r_k))
\geq
c r_k^n.
\end{equation*}
That is, \zcref{lowerdensityradius} holds with $r_k$ replaced by $3r_k$.

It follows by induction that for any \(j\in\mathbb{N}\) with \(3^jr_{k}< \delta R/4\) we have
\[
\mu(B(x,3^jr_{k}))
\geq
c(3^{j-1}r_{k})^n
.
\]
Finally, for any $r<\delta R/4$, since $r_k\to 0$, we may choose \(j,k\) with \(3^jr_{k}\leq r<3^{j+1}r_{k}\) to obtain
\[
\mu(B(x,r))
\geq
c(r/9)^n
.
\]
\end{proof}

We now prove our main result.
\begin{proof}[Proof of \zcref{thm-density}]
Abbreviate \(\mu=\restrict{\calH^n}E\).
Let $X=\bigcup_j X_j$ and $\delta_j,\tau_j>0$ be given by \zcref{filtering}.
Fix $j\in\N$ and let \(\varepsilon_j,c_j>0\) be the parameters from \zcref{pro_inductionstep} corresponding to \(\tau_j,\delta_j\) and $C=3^n$.
Then for any $\varepsilon>0$ there exists a further decomposition $X_j=\cup_k X_j^k$ and for each $k\in\N$ a matrix $P=(e_1,\ldots,e_n)$ of unit vectors with $|P^{-1}|\leq \tau_j$ such that each $\restrict{\mu}{X_j^k}$ has Alberti representations \(\P_1,\ldots,\P_n\) with \(\P_i\in\Gcb{e_i}{\varepsilon}{\delta_j}\),
and by \zcref{phigammabilip}
in particular
\begin{align*}
\restrict \mu{X_j^k}
&\ll
\nFF{\phidiffabs}{\P_i}
.
\end{align*}

We now also fix $k\in\N$.
By \zcref{layer-cake} we can further decompose \(X_j^k=\bigcup_l X_{j}^{k,l}\) such that for each \(l\in\N\) there exist \(\P_i^l\) with
\[
\restrict \mu{X_{j}^{k,l}}
=
\nFF{\phidiffabs}{\P_i^l}
.
\]
Finally fix $l\in\N$.
By the moreover statement of \zcref{layer-cake},
\[
\P_i^l(\Gamma(X)\setminus \Gcb{e_i}{\epsilon}{\delta_j})=0
.
\]
Thus we may apply \zcref{almosteverywheredensity} and conclude $\Theta_{*}^{n}(E\cap X_{j,l}^k,x)>0$ for \(\mu\)-almost every \(x \in X_{j,l}^k\).
\end{proof}

\printbibliography

\end{document}